\title{Neighboring ternary cyclotomic coefficients differ by at most one}
\author{Yves Gallot and Pieter Moree}
\documentclass[12pt]{article}
\usepackage{amssymb, latexsym, amsfonts}
\def\@ptsize{2}
\newtheorem{Thm}{Theorem}

\newtheorem{Lem}{Lemma}

\newtheorem{Def}{Definition}
\newtheorem{cor}{Corollary}

\newcommand{\qed}{\hfill $\Box$}

\begin{document}
\date{}
\maketitle
{\def\thefootnote{}
}
\begin{abstract}
\noindent A cyclotomic polynomial $\Phi_n(x)$ is said to be ternary if $n=pqr$ with $p,q$ and $r$
distinct odd prime factors. Ternary cyclotomic polynomials are the simplest ones for which the
behaviour of the coefficients is not completely understood.
Eli Leher showed in 2007 that neighboring ternary cyclotomic coefficients differ
by at most four. We show that, in fact, they differ by at most one. Consequently, the
set of coefficients occurring in a ternary cyclotomic polynomial consists of consecutive integers.\\
\indent As an application we reprove in a simpler way a result of Bachman from 2004 on ternary
cyclotomic polynomials with an optimally large set of coefficients.
\end{abstract}
\section{Introduction}
The $n$th cyclotomic polynomial $\Phi_n(x)$ is defined by
$$\Phi_n(x)=\prod_{1\le j\le n\atop (j,n)=1}(x-\zeta_n^j)=\sum_{k=0}^{\infty}a_n(k)x^k,$$ with
$\zeta_n$ a $n$th primitive root of unity (one can take $\zeta_n=e^{2\pi i/n}$).
It has degree $\varphi(n)$, with $\varphi$ Euler's totient function. 
We write
$$f(x)=\sum_{k=0}^{\infty}c_kx^k=\sum_{k=0}^{{\rm deg}(f)}c_kx^k,$$
and put ${\cal C}(f)=\{c_k:0\le k\le {\rm deg}(f)\}$ and ${\cal C}_0(f)=
\{c_k:k\ge 0\}$.
Note that ${\cal C}_0(f)={\cal C}(f)\cup \{0\}$. {}For notational convenience we will write
${\cal C}(n)$ instead of ${\cal C}(\Phi_n)$ and ${\cal C}_0(n)$ instead of ${\cal C}_0(\Phi_n)$.
\begin{Def}
If  ${\cal C}_0(n)\subseteq \{-1,0,1\}$, then $\Phi_n$ is said to be flat.
\end{Def}
In the 19th century it was noted that $\Phi_n$ has a strong tendency to be flat and this intrigued various
mathematicians enough to study the coefficients of $\Phi_n$ more intensively. {}For some recent contributions see
e.g. Bachman \cite{B3} and Kaplan \cite{Kaplan}.\\
\indent Using Lemma \ref{binary} below it is not difficult to establish the classical fact that if $\Phi_n$
is not flat, then $n$ has at least three distinct odd prime factors. The simplest case arises when
$n=pqr$ with $2<p<q<r$ odd primes. In this case $n$ is said to be ternary and $\Phi_n$ is said
to be a ternary cyclotomic polynomial. Given $\Phi_{pqr}$ one of the basic problems is to
determine the maximum (in absolute value) of its coefficients. 
Put $M(p)=\max\{|m|:m\in {\cal C}(n),~n=pqr,~p<q<r~{\rm primes~}\}$.
Here it is known \cite{B1} that
$M(p)\le 3p/4$. In 1968 it was conjectured by Sister Marion Beiter \cite{Beiter-1} (see 
also \cite{Beiter-2}) that
$M(p)\le (p+1)/2$. 
She proved it for $p\le 5$.
The first to show that Beiter's conjecture is false seems to have
been Eli Leher \cite[p. 70]{Leher}, who gave the counter-example $a_{17\cdot29\cdot 41}(4801)=-10$,
showing that $M(17)\ge 10>9=(17+1)/2$. The present authors \cite{GM} provided infinitely many 
counter-examples for
the case $p=17$ and in fact for every $p\ge 11$.
Moreover, they have shown that
for every $\epsilon>0$ and $p$ sufficiently large $M(p)>({2\over 3}-\epsilon)p$. Thus only for $p=7$
the Beiter conjecture remains open.\\
\indent In his PhD thesis on numerical semigroups Leher shows
that in case $n$ is ternary one has $|a_n(k)-a_n(k-1)|\le 4$ (\cite[Theorem 57]{Leher}) and remarks 
that he does not know whether this
bound is sharp. Here we show that it is not.
\begin{Thm} 
\label{main}
Let $n$ be ternary, that is $n=pqr$ with $2<p<q<r$ odd primes.
Then, for $k\ge 1$, $|a_n(k)-a_n(k-1)|\le 1$.
\end{Thm}
\begin{cor}
If $n$ is ternary, then ${\cal C}(n)=\{a,a+1,\ldots,b-1,b\}$ with $a$ and $b$ integers,
that is ${\cal C}(n)$ consists of a range of consecutive integers.
\end{cor}
{}For convenience we will say that $f\in \mathbb Z[x]$ has the {\it jump one} property if neighboring
coefficients differ by at most one. Thus Theorem \ref{main} says that a ternary $\Phi_n$ has the jump
one property.
If ${\cal C}_0(n)$ consists of a range of consecutive integers, we say $\Phi_n$ is 
{\it coefficient convex}. Thus Corollary 1 says that a ternary $\Phi_n$ is coefficient convex.
Notice that if $\Phi_n$ is flat, then $\Phi_n$ is coefficient
convex. In Section 6 we consider the problem of determining those $n$ for which
$\Phi_n$ is coefficient convex.\\
\indent Leher uses ideas from the theory of semigroups to prove his result. In Section
2 we discuss the connection between numerical semigroups and cyclotomic polynomials, for
further details we refer to Leher's PhD thesis.\\
\indent In Section 4 we prove Theorem \ref{main}. The proof does not use any semigroup
ideas, but rests on a recent lemma of Kaplan that is discussed in Section 2, along with
some examples.\\
\indent In Section 5 we show how our main result makes detecting so-called optimal ternary
cyclotomic polynomials easier and demonstrate this by giving a reproof of the main result in 
Bachman \cite{B2}. The proof crucially rests on the examples considered in Section 3.1.\\
\indent {}For a nice survey of basic properties of cyclotomic coefficients we refer
to Thangadurai \cite{Thanga}.

\section{Binary cyclotomic polynomials and numerical semigroups}
A number $m$ is said to be a {\it natural combination} of the integers $a_1,\ldots,a_m$ if there are
non-zero integers $k_1,\ldots,k_m$ such that $n=k_1a_1+\cdots+k_ma_m$. Let
$A=\{a_1,\ldots,a_m\}$ be a set of natural numbers and 
$S=S(A)=S(a_1,\ldots,a_m)$ be the set of all natural combinations of the $a_i$'s. Then
$S$ is a {\it semigroup} (that is, it is closed under addition). The semigroup $S$ is
{\it numerical} if its complement $\mathbb Z_{\ge 0}\backslash S$ is finite. If $S$ is numerical,
then $\max \{\mathbb Z_{\ge 0}\backslash S\}=F(S)$ is the {\it Frobenius number} of $S$. 
It is not difficult to prove that $S(a_1,\ldots,a_m)$ is numerical iff $a_1,\ldots,a_m$ 
are relatively prime. The {\it Hilbert series} of the numerical semigroup $S$ is the formal
power series $H_G(x)=\sum_{s\in S}x^s\in \mathbb Z[[x]]$. {}For a numerical semigroup $S$,
$(1-x)H_S(x)$ is a polynomial of degree $F(S)+1$. Leher writes $P_S(x)=(1-x)H_S(x)$ and
calls $P_S(x)$ the {\it semigroup polynomial}. It can be shown that 
$P_{S(p,q)}(x)=\Phi_{pq}(x)$. This leads to the following interpretation of the coefficients
$a_{pq}(k)$:
$$a_{pq}(k)=\cases{1 & if $k\in S(p,q),~k-1\not\in S(p,q)$;\cr
-1 & if $k\not\in S(p,q),~k-1\in S(p,q)$;\cr
0 & otherwise.}$$
Lemma \ref{binary} below gives an explicit evaluation of the binary coefficients $a_{pq}(k)$. 
Leher \cite[Theorem 50]{Leher} shows that an analogous evaluation holds for the coefficients of $P_{S(p,q)}(x)$ in
case $p$ and $q$ are relatively prime positive integers exceeding one.\\
\indent The above material leads to some natural questions. {}For which numerical semigroups
$S$ do we have that $\sum_{x\in S}x^s$ divides $x^m-1$ for some $m$ ? Another question
is to determine all integers $n$ for which
\begin{equation}
\label{elias}
(1-x)\sum_{s\in S_n}x^s=\Phi_n(x),
\end{equation} 
for some set of integers $S_n$.

\section{Kaplan's lemma reconsidered}
\indent Our main tool will be the following recent result due to Kaplan \cite{Kaplan}, the
proof of which uses the identity
$$\Phi_{pqr}(x)=(1+x^{pq}+x^{2pq}+\cdots)(1+x+\cdots+x^{p-1}-x^q-\cdots-x^{q+p-1})
\Phi_{pq}(x^r).$$
\begin{Lem} {\rm (Nathan Kaplan, 2007)}.
\label{kapel}
Let $2<p<q<r$ be primes and $k\ge 0$ be an integer.
Put $$b_i=\cases{a_{pq}(i) & if $ri\le k$; \cr
0 & otherwise.}$$ We have
\begin{equation}
\label{lacheens}
a_{pqr}(k)=\sum_{m=0}^{p-1}(b_{f(m)}-b_{f(m+q)}),
\end{equation}
where $f(m)$ is the unique integer  such that
$f(m)\equiv r^{-1}(k-m) ({\rm mod~}pq)$ and $0\le  f(m) < pq$.
\end{Lem}
This lemma reduces the computation of $a_{pqr}(k)$ to that of $a_{pq}(i)$ for
various $i$. These binary cyclotomic polynomial coefficients are computed
in the following lemma. {}For a proof
see e.g. Lam and Leung \cite{LL} or Thangadurai \cite{Thanga}.
\begin{Lem}
\label{binary}
Let $p<q$ be odd primes. Let $\rho$ and $\sigma$ be the (unique) non-negative
integers for which $1+pq=(\rho+1) p+(\sigma+1) q$.
Let $0\le m<pq$. Then either $m=\alpha_1p+\beta_1q$ or $m=\alpha_1p+\beta_1q-pq$
with $0\le \alpha_1\le q-1$ the unique integer such that $\alpha_1 p\equiv m({\rm mod~}q)$
and $0\le \beta_1\le p-1$ the unique integer such that $\beta_1 q\equiv m({\rm mod~}p)$.
The cyclotomic coefficient $a_{pq}(m)$ equals
$$\cases{1 & if $m=\alpha_1p+\beta_1q$ with $0\le \alpha_1\le \rho,~0\le \beta_1\le
\sigma$;\cr -1 & if $m=\alpha_1p+\beta_1q-pq$ with $\rho+1\le \alpha_1\le q-1,~\sigma+1\le 
\beta_1\le p-1$;\cr  0 & otherwise.}
$$
\end{Lem}
We say that $[m]_p=\alpha_1$ is the {\it $p$-part of $m$} and $[m]_q=\beta_1$ is the {\it $q$-part
of $m$}. It is easy to see that
$$m=\cases{[m]_pp+[m]_qq & if $[m]_p\le \rho$ and $[m]_q\le \sigma$;\cr
[m]_pp+[m]_qq-pq & if $[m]_p>\rho$ and $[m]_q>\sigma$;\cr
[m]_pp+[m]_qq-\delta_mpq & otherwise,}$$
with $\delta_m\in \{0,1\}$. Using this observation we find that, for $i<pq$,
$$b_i=\cases{1 & if $[i]_p\le \rho$, $[i]_q\le \sigma$ and $[i]_pp+[i]_qq\le k/r$;\cr
-1 & if $[i]_p>\rho$, $[i]_q>\sigma$ and $[i]_pp+[i]_qq-pq\le k/r$;\cr
0 & otherwise.}$$
Thus in order to evaluate $a_{pqr}(n)$ using Kaplan's lemma, it is not necessary to compute
$f(m)$ and $f(m+q)$ (as we did in \cite{GM}), it suffices to compute $[f(m)]_p$, $[f(m)]_q$, 
$[f(m+q)]_p$ and $[f(m+q)]_q$ (which is easier). Indeed, as  $[f(m)]_p=[f(m+q)]_p$, it 
suffices to compute $[f(m)]_p$, $[f(m)]_q$, 
and $[f(m+q)]_q$.\\
\indent {}For future reference we provide a version of Kaplan's lemma in which the
computation of $b_i$ has been made explicit, and thus is selfcontained.
\begin{Lem}
\label{kapel2}
Let $2<p<q<r$ be primes and $k\ge 0$ be an integer.
We put $\rho=[(p-1)(q-1)]_p$ and $\sigma=[(p-1)(q-1)]_q$. 
Furthermore, we put $$b_i=\cases{1 & if $[i]_p\le \rho$, $[i]_q\le \sigma$ and $[i]_pp+[i]_qq\le k/r$;\cr
-1 & if $[i]_p>\rho$, $[i]_q>\sigma$ and $[i]_pp+[i]_qq-pq\le k/r$;\cr
0 & otherwise.}$$
We have
\begin{equation}
\label{lacheens2}
a_{pqr}(k)=\sum_{m=0}^{p-1}(b_{f(m)}-b_{f(m+q)}),
\end{equation}
where $f(m)$ is the unique integer  such that
$f(m)\equiv r^{-1}(k-m) ({\rm mod~}pq)$ and $0\le  f(m) < pq$.
\end{Lem}
Note that if $i$ and $j$ have the same $p$-part, then $b_ib_j\ne -1$, that is $b_i$ and $b_j$
cannot be of opposite sign. {}From this it follows that $|b_{f(m)}-b_{f(m+q)}|\le 1$, and
thus we infer from Kaplan's lemma that $|a_{pqr}(k)|\le p$. It is possbile to improve
on this argument and get a sharper bound for $|a_{pqr}(k)|$. We hope to return to this issue in 
a future paper. Of course if $i$ and $j$ have the same $q$-part, then $b_ib_j\ne -1$ also.
\subsection{Examples of computing coefficients with Kaplan's lemma}
In this section we carry out a sample computation using Kaplan's lemma. {}For more involved examples the
reader is referred to \cite{GM}.\\
\indent We remark that the result that $a_n(k)=(p+1)/2$ in Lemma \ref{moellerext} is due
to Herbert M\"oller \cite{HM}. The reproof we give is rather different. The foundation for M\"oller's result
is due to Emma Lehmer \cite{Emma}, who already in 1936 had shown that $a_n({1\over 2}(p-3)(qr+1))=(p-1)/2$ with
$p,q,r$ and $n$ satisfying the conditions of Lemma \ref{moellerext}.
\begin{Lem}
\label{moellerext}
Let $p<q<r$ be primes satisfying
$$p>3,~q\equiv 2({\rm mod~}p),~r\equiv {p-1\over 2}({\rm mod~}p),~r\equiv {q-1\over 2}({\rm mod~}q).$$
Put $n=pqr$ and $k=(p-1)(qr+1)/2$. Then $a_{n}(k-r)=-(p-1)/2$ and $a_n(k)=(p+1)/2$.
\end{Lem}
{\it Proof}. First it will be shown that $a_n(k)=(p+1)/2$.
Using that $q\equiv 2({\rm mod~}p)$, we infer from
$1+pq=(\rho+1)p+(\sigma+1)q$ that $\sigma={p-1\over 2}$ and $(\rho+1)p=1+({p-1\over 2})q$ (and
hence $\rho=(p-1)(q-2)/(2p)$). On invoking the Chinese remainder theorem on checks that
\begin{equation}
\label{LD}
-{1\over r}\equiv 2\equiv -({q-2\over p})p+q({\rm mod~}pq).
\end{equation}
Furthermore, writing $f(0)$ as a linear combination of $p$ and $q$ we see that
\begin{equation}
\label{efnul}
f(0)\equiv {k\over r}\equiv ({p-1\over 2})q+{p-1\over 2r}\equiv ({p-1\over 2})q+1-p\equiv \rho p({\rm
mod~}pq).
\end{equation}
{}From (\ref{LD}) and (\ref{efnul}) we infer that, for $0\le m\le (p-1)/2$, we have
$[f(m)]_p=\rho-m(q-2)/p\le \rho$ and $[f(m)]_q=m\le \sigma$. On noting that
$[f(m)]_pp+[f(m)]_qq=\rho p+2m\le \rho p+ p-1=[k/r]$, we infer that
$a_{pq}(f(m))=b_{f(m)}=1$ in this range (see also Table 1).\\

\centerline{\bf TABLE 1}
\begin{center}
\begin{tabular}{|c|c|c|c|c|c|c|}
\hline
$m$ & $[f(m)]_p$ & $[f(m)]_q$ & $f(m)$ & $a_{pq}(f(m))$ & $b_{f(m)}$\\
\hline
0 & $\rho$ & $0$ & $\rho p$ & 1 & 1\\
\hline
1 & $\rho-(q-2)/p$ & 1 & $\rho p+2$ & 1 & 1\\
\hline
$\ldots$ & $\ldots$ & $\ldots$ & $\ldots$ & 1 & 1\\
\hline
$j$ & $\rho-j(q-2)/p$ & $j$ & $\rho p +2j$ & 1 & 1\\
\hline
$\ldots$ & $\ldots$ & $\ldots$ & $\ldots$ & 1 & 1\\
\hline
$(p-1)/2$ & $0$ & $(p-1)/2$ & $(p-1)q/2$ & 1 & 1\\
\hline
\end{tabular}
\end{center}

Note that $f(m)\equiv f(0)-m/r\equiv \rho p+2m({\rm mod~}pq)$, from which one easily
infers that $f(m)=\rho p+2m$ for $0\le m\le p-1$ (as $\rho p+2m\le \rho p+2(p-1)<pq$). In
the range ${p+1\over 2}\le m\le p-1$ we have $f(m)\ge \rho p+p+1=(p-1)q/2+2>k/r$, and hence
$b_{f(m)}=0$.\\
\indent On noting that $f(m+q)\equiv f(m)-q/r\equiv f(m)+2q\equiv \rho p+2m+2q({\rm mod~}pq)$,
one easily finds, for $0\le m\le p-1$, that $f(m+q)=\rho p+2m+2q>k/r$ and hence $b_{f(m+q)}=0$.\\
\indent By Kaplan's lemma one infers that
$$a_n(k)=\sum_{m=0}^{p-1}(b_{f(m)}-b_{f(m+q)})=\sum_{m=0}^{(p-1)/2}1={p+1\over
2}.$$
\indent Next we show that $a_n(k-r)=-(p-1)/2$. Put $k'=k-r$ and
$f_1(m)\equiv r^{-1}(k-r-m)({\rm mod~}pq)$ with $0\le f_1(m)<pq$. Furthermore we put
$b''_{f_1(m)}=a_{pq}(f_1(m))$ if $f_1(m)\le k'/r$ and zero otherwise. By (\ref{efnul}) we deduce
that
$$f_1(0)\equiv f(0)-1\equiv (q-1)p+({p-1\over 2})q({\rm mod~}pq).$$
An easy calculation yields the correctness of Table 2.\\

\centerline{\bf TABLE 2}
\begin{center}
\begin{tabular}{|c|c|c|c|c|c|}
\hline
$m$ & $[f_1(m)]_p$ & $[f_1(m)]_q$ & $f_1(m)$ & $b''_{f_1(m)}$\\
\hline
0 & $q-1$ & $(p-1)/2$ & $(p-1)q/2-p$ & 0\\
\hline
1 & $q-1-(q-2)/p$ & $(p+1)/2$ &  $(p-1)q/2-p+2$ & -1\\
\hline
$\ldots$ & $\ldots$ & $\ldots$ & $\ldots$ & -1\\
\hline
$j$ & $q-1-j(q-2)/p$ & $(p-1)/2+j$ & $(p-1)q/2-p+2j$ & -1\\
\hline
$\ldots$ & $\ldots$ & $\ldots$ & $\ldots$ & -1\\
\hline
$(p-1)/2$ & $q-1-\rho$ & $p-1$ & $(p-1)q/2-1$ & -1\\
\hline
$(p+1)/2$ & $q-1-\rho-(q-2)/p$ & 0 & $(p-1)q/2+1$ & 0\\
\hline
\end{tabular}
\end{center}

{}For $(p+1)/2\le m\le p-1$ one finds that $f_1(m)>k'/r$ and hence $b''_{f_1(m)}=0$.
On noting that, for $0\le m\le p-1$, 
$f_1(m+q)=({p-1\over 2})q-p+2m+2q>k'/r$, we find that $b''_{f_1(m+q)}=0$ in this 
range. Kaplan's lemma now gives $$a_n(k-r)=\sum_{m=0}^{p-1}(b''_{f_1(m)}-
b''_{f_1(m+q)})=-\sum_{m=1}^{(p-1)/2}1=-{(p-1)\over 2},$$
completing the proof. \qed\\

\noindent We recall from \cite{GM} the following result, which is implicit in
Kaplan's paper \cite{Kaplan}.
\begin{Lem} Let $2<p<q<r$ be primes and $n\ge 0$ be an integer.
Suppose that $a_{pqr}(n)=m$. Write
$n=[{n\over r}]r+n_0$ with $0\le n_0<r$. 
Let $t>pq$ be a prime satisfying $t\equiv -r({\rm mod~}pq)$. Let
$0\le n_1<pq$ be the unique integer such that 
$n_1\equiv q+p-1-n_0({\rm mod~}pq)$. Then
$$a_{pqt}\Big(\Big[{n\over r}\Big]t+n_1\Big)=-m.$$
\end{Lem}
Using the latter lemma one immediately gets from Lemma \ref{moellerext} the following one, however
with the condition $r_1>q$ replaced by $r_1>pq$. 
On proceeding as in the proof of
Lemma \ref{moellerext}, one gets the feeling of d\'ej\`a vu. Indeed, it turns out that
on running through $m=0,\ldots,p-1$, $m=q,q+1,\ldots,q+p-1$, the $f(m)$ in the setup 
of Lemma \ref{moellerext} correspond to the $f(m)$ for $q+p-1-m$ in the setup of
Lemma \ref{bloeh}, the effect being that the corresponding ternary coefficients differ
by a minus sign. On doing this it turns out that the condition $r_1>pq$ can be relaxed
to the condition $r_1>q$.
\begin{Lem}
\label{bloeh}
Let $p<q<r$ be primes satisfying
$$p>3,~q\equiv 2({\rm mod~}p),~r_1\equiv {p+1\over 2}({\rm mod~}p),~r_1\equiv {q+1\over 2}({\rm mod~}q).$$
Put $n_1=pqr_1$ and $k_1=(p-1)(qr_1+1)/2+q$. Then $a_{n_1}(k_1-r_1)=(p-1)/2$ and $a_{n_1}(k_1)=-(p+1)/2$.
\end{Lem}
{\tt Remark}. Note that, in Lemma \ref{moellerext}, $r=(lpq-1)/2$ for some odd $l\ge 1$ and
that $r_1=(l_1pq+1)/2$ for some odd $l_1\ge 1$.

\section{Proof of the jump one property}
Having Kaplan's lemma at our disposal we are ready to give a proof of the jump one property.\\
{\it Proof of Theorem \ref{main}}. Let $f'(m)$ be the unique integer $0\le  f'(m) < pq$ such that
$f'(m)\equiv r^{-1}(k-1-m) ({\rm mod~}pq)$. Let $b'_i$ be defined as $b_i$, but 
with $k$ replaced by $k-1$. 
Note that
\begin{equation}
\label{divvie}
b_i-b'_i=\cases{b_{k\over r} & if $r|k$ and $i={k\over r}$;\cr
0 & otherwise.}
\end{equation}
If $r|k$, then $f(0)=k/r$. {}For $j=0,\ldots,p-2$ we have
$f'(j)=f(j+1)$ and since $f(j+1)\ne f(0)$, we infer using (\ref{divvie}) that
$$\gamma_1:=\sum_{j=0}^{p-1}\Big(b_{f(j)}-b'_{f'(j)}\Big)=b_{f(0)}-b'_{f'(p-1)}.$$
Likewise we see that
$$\gamma_2:=\sum_{j=0}^{p-1}\Big(b_{f(j+q)}-b'_{f'(j+q)}\Big)=b_{f(q)}-b'_{f'(p-1+q)}.$$
By Kaplan's lemma it then follows that
\begin{equation}
\label{lulu}
a_{n}(k)-a_{n}(k-1)=\gamma_1-\gamma_2=b_{f(0)}-b_{f(q)}-b'_{f'(p-1)}+b'_{f'(p-1+q)}.
\end{equation}
Denote $f(0),f(q),f'(p-1),f'(p-1+q)$ by, respectively, $\alpha_1,\alpha_2,\alpha_3,\alpha_4$.
Note that modulo $pq$ we have
$$\alpha_1\equiv {k\over r},~\alpha_2\equiv {k-q\over r},~\alpha_3\equiv 
{k-p\over r},~\alpha_4\equiv {k-p-q\over r}.$$
Denote $b_{f(0)}, b_{f(q)}, b'_{f'(p-1)}$, $b'_{f'(p-1+q)}$ by, respectively, 
$\beta_1,\beta_2,\beta_3$ and $\beta_4$. Thus we can rewrite (\ref{lulu}) as
\begin{equation}
\label{lala}
a_{n}(k)-a_{n}(k-1)=\beta_1-\beta_2-\beta_3+\beta_4.
\end{equation}
Since 
$\alpha_1$ and $\alpha_2$ have equal $p$-part, $\beta_1$ and $\beta_2$ cannot
be of opposite sign and hence $|\beta_1-\beta_2|\le 1$ and by a similar argument
we find $|\beta_3-\beta_4|\le 1$. It follows that $|a_{n}(k)-a_{n}(k-1)|\le 2$.\\
\indent Since $\alpha_1$ and $\alpha_3$ have the same $q$-part, $\beta_1$ and
$\beta_3$ cannot be of opposite
sign. Likewise, $\beta_2$ and $\beta_4$ cannot be of opposite sign. Put ${\overline
\beta}=(\beta_1,\beta_2,\beta_3,\beta_4)$.
It follows that if
$a_{n}(k)-a_{n}(k-1)=2$, then 
\begin{equation}
\label{bie}
{\overline \beta}=(1,0,0,1)  {\rm ~or~} {\overline \beta}=(0,-1,-1,0). 
\end{equation}
Likewise we infer that if
$a_{n}(k)-a_{n}(k-1)=-2$, then 
\begin{equation}
\label{boe}
{\overline \beta}=(-1,0,0,-1)  {\rm ~or~} {\overline \beta}=(0,1,1,0). 
\end{equation}
The proof is completed if we can show that none of these four possibilities for ${\overline \beta}$
can occur.\\
-{\tt Excluding}  ${\overline \beta}=(1,0,0,1)$: We have
$\alpha_1=[\alpha_1]_pp+[\alpha_1]_qq$.and $\alpha_4=[\alpha_4]_pp+[\alpha_4]_qq$ with
$[\alpha_1]_p\le \rho, [\alpha_1]_q\le \sigma, [\alpha_4]_p\le \rho, [\alpha_4]_q\le \sigma$, 
$\alpha_1\le k/r$ and $\alpha_4\le (k-1)/r$. Since $\alpha_2$ has the same $p$-part as
$\alpha_1$ and the same $q$-part as $\alpha_4$, we find that
$\alpha_2\equiv [\alpha_1]_pp+[\alpha_4]_qq({\rm mod~}pq)$. Since
$[\alpha_1]_pp+[\alpha_4]_qq\le \rho p+\sigma q<pq$, it follows that 
$\alpha_2=[\alpha_1]_pp+[\alpha_4]_qq$. Since by assumption $\beta_2=0$, we must
have $\alpha_2>k/r$ (if $\alpha_2\le k/r$, then we would have $\beta_2=1$). Likewise we 
infer that $\alpha_3=[\alpha_4]_pp+[\alpha_1]_qq$ and
$\alpha_3>(k-1)/r$. It follows that $\alpha_2+\alpha_3>(2k-1)/r$. On the other hand,
$\alpha_2+\alpha_3=\alpha_1+\alpha_4\le (2k-1)/r$. This contradiction shows that the
case ${\overline \beta}=(1,0,0,1)$ cannot occur.\\
-{\tt Excluding}  ${\overline \beta}=(0,-1,-1,0)$: We have 
$$\alpha_2=[\alpha_2]_pp+[\alpha_2]_qq-pq,~\alpha_3=[\alpha_3]_pp+[\alpha_3]_qq-pq$$ with
$[\alpha_2]_p>\rho, [\alpha_2]_q> \sigma, [\alpha_3]_p>\rho, [\alpha_4]_q >\sigma$, 
$\alpha_2\le k/r$ and $\alpha_3\le (k-1)/r$. Since $\alpha_1$ has the same $p$-part as
$\alpha_2$ and the same $q$-part as $\alpha_3$, we find that
$\alpha_1\equiv [\alpha_2]_pp+[\alpha_3]_qq-pq({\rm mod~}pq)$. Since
$0\le [\alpha_2]_pp+[\alpha_3]_qq-pq<pq$, we infer that 
$\alpha_1=[\alpha_2]_pp+[\alpha_3]_qq-pq$. Since by assumption $\beta_1=0$, we must
have $\alpha_1>k/r$ (for otherwise we would have $\beta_1=-1$). Likewise we infer that
$\alpha_4=[\alpha_3]_pp+[\alpha_2]_qq-pq$. On the one hand we have
$\alpha_2+\alpha_3\le (2k-1)/r$, on the other hand we have $\alpha_2+\alpha_3=\alpha_1+\alpha_4>(2k-1)/r$,
a contradiction showing that ${\overline \beta}=(0,-1,-1,0)$ cannot occur.\\
-{\tt Excluding the two remaining cases}: can be done by minor variations of the above arguments and is
left to the interested reader.\\
\indent Thus the proof is completed. \qed\\

\noindent {\tt Remark}. The result of Leher that $|a_n(k)-a_n(k-1)|\le 4$ is of course an immediate
consequence of  (\ref{lulu}).

\section{Coefficient optimal ternary polynomials}
In this section we give an application of the jump one property.\\
\indent The difference between the largest and the smallest coefficients of $\Phi_{pqr}$ is known
to be at most $p$ \cite[(1.5)]{B2}. We say that $\Phi_{pqr}$ is coefficient optimal if the difference
between the largest and smallest coefficient is exactly $p$. Bachman \cite{B2} found two infinite
families of coefficient optimal ternary polynomials $\Phi_{pqr}$, with
${\cal C}(pqr)=[-(p-1)/2,(p+1)/2]$ 
for one family and ${\cal C}(pqr)=[-(p+1)/2,(p-1)/2]$
for the other family. Using the jump one property one immediately infers the following
result.
\begin{Lem} 
\label{vijf}
Let $p<q<r$ be odd primes.
If $a,b\in {\cal C}(pqr)$ and $b-a=p$, then ${\cal C}(pqr)=\{-a,a+1,\ldots,b-1,b\}$.
\end{Lem}
Thus the jump one property might be helpful in studying families of coefficient optimal
ternary cyclotomic polynomials. We demonstrate this by showing how it can be used to
reprove the main result in Bachman (whose proof is very different).
\begin{Thm}
Let $p,q,r,k,n$ be as in 
Lemma {\rm \ref{moellerext}}. Then $\Phi_n$ is coefficient optimal and in particular
\begin{equation}
\label{aa}
{\cal C}(n)=\{a_{n}(k-r),a_{n}(k-r+1),\ldots,a_{n}(k)\}=[-(p-1)/2,(p+1)/2]\cap \mathbb Z.
\end{equation}
Let $p,q,r_1,k_1,n_1$ be as in 
Lemma {\rm \ref{bloeh}}.
Then $\Phi_{n_1}$ is coefficient optimal and in particular
\begin{equation}
\label{ab}
{\cal C}(n_1)=\{a_{n_1}(k_1-r_1),\ldots,a_{n_1}(k_1)\}=[-(p+1)/2,(p-1)/2]\cap \mathbb Z.
\end{equation}
\end{Thm}
{\it Proof}. By Lemma \ref{moellerext} one has $a_n(k-r)=-(p-1)/2$ and
$a_n(k)=(p+1)/2$. By the jump one property it then follows that the second equality
in (\ref{aa}) holds. By Lemma \ref{vijf} with $b=(p+1)/2$ and $a=-(p-1)/2$ it follows
that the first equality in (\ref{aa}) holds.\\
\indent The proof of the remaining assertion (\ref{ab}) is completely similar, but makes
use of Lemma \ref{bloeh} instead of Lemma \ref{moellerext}. \qed

\section{Coefficient convexity}
Theorem \ref{main} gives naturally rise to the notion of (strong) coefficient convexity. In
this section we consider the issue of coefficient convexity of cyclotomic(-like) polynomials
in somewhat greater detail, leaving the proofs for a future publication \cite{DM}.
\begin{Def}
We say that $\Phi_n$ is coefficient convex if ${\cal C}_0(n)=I_n\cap \mathbb Z$ for some
interval $I_n$ in the reals. We say it is strongly coefficient convex 
if ${\cal C}(n)=I_n\cap \mathbb Z$ for some
interval $I_n$ in the reals.
\end{Def}
A cyclotomic polynomial can be coefficient convex without being strongly coefficient convex, e.g.
$\Phi_{2p}$ ($p$ being an odd prime) is coefficient convex but not strongly coefficient convex. The latter cyclotomic polynomial
also shows that a cyclotomic polynomial can be flat without being strongly coefficient convex.
Moreover, Theorems \ref{bloop} and \ref{cloop} below are false if one replaces `coefficient convex'
by `strongly coefficient convex'.\\
\indent Using Theorem \ref{main} it is not difficult to establish the following result.
\begin{Thm} 
\label{bloop}
Suppose that $n$ has at most $3$ prime factors, then $\Phi_n$ is
coefficient convex.
\end{Thm}
Numerical computations suggest that if $\Phi_n$ is ternary, then $\Phi_{2n}$ is coefficient
convex. If this would be true, then in Theorem \ref{bloop} one can replace `3 prime factors'
by `3 distinct odd prime factors'. This is best possible as the following examples show:\\
$n=7735=5\cdot 7\cdot 13\cdot 17$, ${\cal C}(n)=[-7,5]-\{-9\}$\\
$n=530689=17\cdot 19\cdot 31\cdot 53$, ${\cal C}(n)=[-50,52]-\{-48,47,48,49,50,51\}$.\\
(Here we write $[-a,b]$ for the range of integers $[-a,b]\cap \mathbb Z$.)\\
\indent We note that if $n$ is ternary, then $\Phi_{2n}$ in general does not have the jump
one property.
\begin{Lem}
If the ternary polynomial $\Phi_n$ is not flat, then $\Phi_{2n}$ does not have the jump
one property.
\end{Lem}
{\it Proof}. Suppose that $a_n(k)=m$ and $|m|>1$. Then by Theorem \ref{main} and
the identity $\Phi_{2n}(x)=\Phi_n(-x)$, we infer that
$$|a_{2n}(k)-a_{2n}(k-1)|=|a_n(k)+a_n(k-1)|\ge 2|m|-1>1,$$
completing the proof. \qed\\
\indent Put $$\Psi_n(x)={x^n-1\over \Phi_n(x)}=\prod_{1\le j\le n\atop (j,n)>1}(x-\zeta_n^j).$$
Write $\Psi_n(x)=\sum_{k=0}^{n-\varphi(n)}c_n(k)x^k$. The coefficients $c_n(k)$ are
integers that turn out to behave in a way quite similar to the cyclotomic coefficients
$a_n(k)$. Apparently Moree \cite{Moree} was the first to systematically study these
coefficients, which he called inverse cyclotomic polynomial coefficients. Here it is
not difficult to prove the following result.
\begin{Thm} 
\label{cloop}
Suppose that $n$ has at most $3$ distinct odd prime factors, then $\Psi_n$ is
coefficient convex.
\end{Thm}
If $n$ has four or more distinct odd prime factors, then $\Psi_n$ need not be
coefficient convex, since we have for example\\
$n=60095=5\cdot 7\cdot 17\cdot 101$, ${\cal C}(\Psi_n)=[-12,12]-\{-11,11\}$.\\
$n=207805=5\cdot 13\cdot 23\cdot 139$, ${\cal C}(\Psi_n)=[-16,16]-\{-15,-13,13,15\}$.\\
$n=335257=13\cdot 17\cdot 37\cdot 41$, ${\cal C}(\Psi_n)=[-40,40]-\{-39,-37,-36,37,39\}$.\\
\indent As we have seen the polynomials $\Phi_n(x)$ 
and $\Psi_n(x)$ are divisors of $x^n-1$ that have the tendency to be
coefficient convex. One can wonder to what extent other divisors
of $x^n-1$ have the same tendency. (Notice that any divisor of $x^n-1$ can be written
as a product of cyclotomic polynomials.)
This problem will be considered in \cite{DM}.

\medskip\noindent {\footnotesize 12 bis rue Perrey,\\
31400 Toulouse, France.\\
e-mail: {\tt galloty@orange.fr}}\\

\medskip\noindent {\footnotesize Max-Planck-Institut f\"ur Mathematik,\\
Vivatsgasse 7, D-53111 Bonn, Germany.\\
e-mail: {\tt moree@mpim-bonn.mpg.de}}
\vskip 5mm

\end{document}